\documentclass[12pt]{article}
\usepackage{amsfonts,amsmath,bbm}
\usepackage[francais]{babel}
\usepackage[latin1]{inputenc}
\usepackage[T1]{fontenc}
\newtheorem{theorem}{Th\'eor\`eme}
\newtheorem{corollary}[theorem]{Corollaire}

\newtheorem{lemma}[theorem]{Lemme}
\newtheorem{proposition}[theorem]{Proposition}
\newtheorem{definition}[theorem]{D\'efinition}

\newtheorem{remark}{Remarque}

\newenvironment{proof}[1][]
    {\par\medbreak{\noindent\bfseries Proof#1.\quad}}
    {\hbox{}\hfill$\square$\bigbreak}

\def\N{\mathbb{N}}
\def\R{\mathbb{R}}
\def\P{\mathcal{P}}

\def\sfleche{\mathcal{S}^{\downarrow}}

\def\eqlaw{\overset{\mathcal{L}}{=}}

\newcommand{\convtenzero}{\underset{t \searrow 0}{\rightarrow}}
\def\ft{\mathcal{F}_t}

\begin{document}

\title{Ranked Fragmentations}
\author{Julien Berestycki}
\date{}

%

\maketitle
\medskip

\noindent {\small Universit\'e de Provence, LATP UMR 6632, CMI, 39
rue F. Joliot-Curie 13453 Marseille cedex 13, France. Email :
jberest@cmi.univ-mrs.fr }

\begin{abstract}
In this paper we define and study self-similar ranked
fragmentations. We first show that any ranked fragmentation is the
image of some partition-valued fragmentation, and that there is in
fact a one-to-one correspondence between the laws of these two types
of fragmentations. We then give an explicit construction of
homogeneous ranked fragmentations in terms of Poisson point
processes. Finally we use this construction and classical results on
records of Poisson point processes to study the small-time behavior
of a ranked fragmentation.
\end{abstract}

\section{Introduction}

Splitting models are meant to describe an object that falls apart.
Applications are numerous and may be found in various fields such
as physical chemistry (aerosols, phase separation,
polymerization), mathematical population genetics or astronomy (we
refer to \cite{aldous-survey} for a survey on applications and
motivations).

This paper focuses on self-similar ranked fragmentation. For the
sake of describing our results, let us just give some heuristic
descriptions while precise definitions will be given in the next
sections.

Imagine a unit-mass object that fragments as time runs. We only
consider the ordered sequence of the fragments masses of this
object so the state space is $$ \sfleche := \{s=(s_1,s_2,...) ,
s_1 \geq s_2 \geq ... \geq 0 , \sum_i s_i \leq 1 \}, $$ the
situation where $\sum_i s_i <1$ corresponding to the fact that a
part of the initial mass has been lost, i.e. the sum of the masses
of the remaining fragments is less than the original total mass.

Let $\lambda=(\lambda(t),t\geq0)$ be a Markov process with values
in $\sfleche$. Call $\lambda$ a self-similar ranked fragmentation
if it fulfills the \textit{scaling} and the \textit{fragmentation}
property.

The \textit{scaling} property means that there exists a real
number $\alpha$, called the index of self-similarity, such that if
$\mathbb{P}_r$ is the law of $\lambda$ started from $(r,0,0,..)$
then the distribution of $(r\lambda(r^{\alpha}t),t\geq 0)$ under
$\mathbb{P}_1$ is $\mathbb{P}_r$.

The \textit{fragmentation} property is a version of the branching
property i.e. for any $u,t\geq 0$ , for any $s= (s_1,s_2,..) \in
\sfleche$, conditionally on $\lambda(u)=s$, $\lambda(t+u)$ has the
same distribution as the variable obtained by concataining and
ordering the sequences $\lambda^{(1)},\lambda^{(2)},..$ where for
each $i$, $\lambda^{(i)}$ has the distribution of $\lambda(t)$
under $\mathbb{P}_{s_i}$.

Here is a simple prototype taken from Brennan and Durrett
\cite{Brennan-Durrett,Brennan-Durrett1} who consider the following
model for polymer degradation : A particle of mass $m$ splits with
exponential rate $m^{\alpha}$, $\alpha \in \R^+$, and gives rise
to two particles of mass $Vm$ and $(1-V)m$, where $V$ is a random
variable with values in $(0,1)$ independent of the past. The new
particles follow the same dynamic independently. The ordered
sequence of the particles masses is a self-similar ranked
fragmentation of index $\alpha$.

This example can be extended in two ways. First one can suppose
that when a particle splits, it might give birth to any number of
particles, possibly infinite, and not just two. Second, in the
example of Brennan and Durett, the splitting times are "discrete",
the first time of splitting is almost surely strictly positive. It
is natural to consider more generally the case where fragmentation
may occur continuously. For instance this happens for the
fragmentation process obtained by logging the continuous random
tree of Aldous in \cite{std-add-coal}.

In the existing literature, a somewhat different class of
processes has been considered, the so-called \textit{partition
valued fragmentation}. Roughly speaking a partition fragmentation,
say $\Pi(t)$, is a process that lives in the space of partitions
of $\N$, such that for any $0<s \leq t$, $\Pi(t)$ is a refinement
of $\Pi(s)$. A way to construct such a fragmentation which makes
clear the connection with the above particle model is the
following : imagine an object $E$ endowed with a unit mass measure
$\mu$ that falls apart as time runs, call \textit{object
fragmentation} the process $F(t)$ with values in partitions of $E$
that describes this fragmentation. Next, let $(u_i)_{i \in \N}$ be
a sequence of iid $E$-valued variables with distribution $\mu$ and
for each $t$ let $\Pi_F(t)$ be the partition of $\N$ such that for
all $i$ and $j$ in $\N$, $i$ and $j$ belong to the same block of
$\Pi(t)$ iff $u_i$ and $u_j$ are in the same fragment of $E$ at
time $t$. By the SLLN we can recover the mass of a fragment as the
asymptotic frequency of the corresponding block. Then $\Pi_F$ is a
partition fragmentation.

Using partition fragmentations to construct ranked fragmentations
is typical of the existing results. These constructions benefit
from two important features :  there is a clear genealogical
structure, and partition fragmentations are characterized by an
index of self-similarity $\alpha$ and a so-called characteristic
exchangeable measure, on which results concerning exchangeability
can be usefully applied (see \cite{aldous} for a survey on
exchangeability).

However partition-valued fragmentations are perhaps less natural
and could be less general than ranked fragmentations, precisely
because we have endowed it with this extra genealogic structure.
In other words it is not clear that an arbitrary ranked
fragmentation can be studied through partition fragmentations.

In section 2 we show that it is in fact the case, and more
precisely that for any ranked fragmentation $\lambda$ we can
associate a partition fragmentation $\Pi$ such that the asymptotic
frequencies of $\Pi$ has same distribution as $\lambda$.

In the next section we use this equivalence between ranked and
partition fragmentations to give a Poisson construction of
homogeneous ranked fragmentation which is an analogue of that
given in \cite{bertoin1} for partition fragmentations. The
difficulty comes from the fact that we can no longer use a
genealogic structure, which played a crucial role in the
partition case.

In section 4, this construction allows us to tackle the study of
small time behavior of a ranked fragmentation. We show that the
2nd largest fragment, correctly renormalized, behaves as the
record of the size of the particles detaching from the main
fragment.

\section{Definitions and first properties}

\subsection{Ranked Fragmentations}

For each $l$ in $[0,1]$ let $P(l)$ be a probability on $\{s \in
\sfleche : \sum_i s_i \leq l\}$ the space of all the possible
fragmentations of $l$. Then for $L=(l_1,l_2,...)$, define $P(L)$
as the distribution on $\sfleche$ of the concatenation and the
decreasing rearrangement of independent $\sfleche$-valued
variables with respective law $P(l_i)$. Call $(P(L) , L \in
\sfleche)$ a fragmentation kernel on $\sfleche$. One says that
the family $(P(l),l\in[0,1])$ generates $(P(L),L \in \sfleche)$.

\begin{definition} \label{def1}
An $\sfleche$-valued process $\lambda(.) $ is called a
$\sfleche$-fragmentation if it is a time-homogeneous Markov
process such that
\begin{enumerate}
\item $\lambda$ is continuous in probability and starts from
$\lambda(0)=(1,0,0,...)$ a.s.
\item the transition semigroup $(P_t(L))$ of $\lambda$ is given by
fragmentation kernels.

\end{enumerate}

\end{definition}

In words, at a given time $t$, each fragment of
$\lambda(t)=(\lambda_1(t),\lambda_2(t),...)$, say $\lambda_i(t)$,
gives rise to an independent fragmentation process which
distribution only depends on the value $\lambda_i(t)$. $\lambda$
is the concatenation and the reordering of all those processes.

For $l \in [0,1]$, let $g_l$ be the application from $\sfleche
\rightarrow \sfleche$ defined by
$$g_l : x=(x_1,x_2,..) \rightarrow (lx_1,lx_2,...).$$

\begin{definition}
\label{Self-similar ranked fragmentation} The fragmentation
$\lambda$, with transition kernels generated by the family
$(P_t(l) ; t \geq 0 , l \in [0,1] )$ is said to be self-similar
with index $\alpha \in \R$ if (in the notations introduced above)
for all $l \in [0,1]$ the distribution $P_t(l)$ coincides with the
image of $P_{l^{\alpha}t}(1)$ by $g_l$.

When $\alpha=0$ the fragmentation is said to be homogeneous.
\end{definition}

$\sfleche$ is endowed with the uniform distance. Note that for any
$s=(s_1,s_2,...) \in \sfleche$ we must have for every $ k \in \N,
s_k \leq \frac{1}{k},$ and thus the uniform and pointwise
convergences are the same. In this setting we prove that a
self-similar $\sfleche$-fragmentation has the Feller property.

\begin{proposition} \label{feller}
(Feller property) The semi-group $P_t$ of a self-similar ranked
fragmentation of index $\alpha$, fulfills the Feller property.
That is $\forall t \geq 0 $ the map
$$L \rightarrow P_t(L)$$ is continuous on $\sfleche$ and for each
fixed $L \in \sfleche$, $P_t(L)$ converge to the Dirac mass at $L$
as $t\rightarrow 0$.
\end{proposition}
\begin{proof}

Consider a sequence $(L_n, n \in \N)$ in $\sfleche$ which
converges to $L \in \sfleche$. Note
$L_n=(l_1^{(n)},l_2^{(n)},...)$, then for all $k$, $l_k^{(n)}
\rightarrow l_k$.

Let $(Y_i(t))_{i \in \N}$ be a sequence of iid
$\sfleche$-fragmentation with same semi-group $(P_{S}(t), S \in
\sfleche ,t \geq 0)$, then, by definition, for all $n \in \N$ the
$\sfleche$ random variable $Z^{(n)}(t)$, obtained by the
decreasing rearrangement of the terms
$g_{l_i^{(n)}}\left(Y_1\left(t(l_i^{(n)})^{\alpha}\right)\right)$
for $i$ in $\N$ :
$$ Z^{(n)}(t)= \left(
g_{l_1^{(n)}}\left(Y_1\left(t(l_1^{(n)})^{\alpha}\right)\right),g_{l_2^{(n)}}\left(Y_2\left(t(l_2^{(n)})^{\alpha}\right)\right),...\right)^{\downarrow}$$
has law $P_{L^{(n)}}(t)$. In the same way $$Z(t)= \left(
g_{l_1}(Y_1(t(l_1)^{\alpha})),g_{l_2}(Y_2(t(l_2)^{\alpha})),...\right)^{\downarrow}$$
has law $P_{L}(t)$.  Now fix $\epsilon>0$, and take
$N>\frac{2}{\epsilon}$. Then $$\forall k\geq N , \forall n \in \N ,
l^{(n)}_k < \epsilon/2.$$ Thus for all $\omega$ $$\sup_{k \geq N }
\left(
\text{dist}\left(g_{l_k^{(n)}}(Y_k(t(l_k^{(n)})^{\alpha})),g_{l_k}(Y_k
(t(l_k)^{\alpha}))\right) \right)<\epsilon.$$ On the other hand, by
the continuity in probability of the processes \\ $(Y_i)_{i\in
\{1,\dots,N-1\}}$, we have that almost surely $$P\left[ \sup_{k \in
\{1,...,N-1\}} \left(
\text{dist}\left(g_{l_k^{(n)}}(Y_k(t(l_k^{(n)})^{\alpha})),g_{l_k}(Y_k
(t(l_k)^{\alpha}))\right) \right)>\epsilon \right]
\underset{n\rightarrow \infty}{\rightarrow} 0.$$ Thus almost surely,
for all $\epsilon
>0$, there exists $N \in \N$ such that for all $n \geq N$
$$\text{dist}\left(Z^{(n)}(t),Z(t) \right)<\epsilon.$$ There is
convergence in probability and thus in law.
\end{proof}

\subsection{Partition Fragmentations}

Most of the results on fragmentation available in the literature
are (or can be) formulated in term of a type of fragmentation
called \textit{partition fragmentation}, which is basically a
process which can be described as a partition of $\N$ getting
finer as time runs.

More precisely, call a subset of $\N$, say $B$, a "block". When
the limit $$|B|:=\lim_{n \rightarrow \infty} \frac{1}{n}
\text{Card} \{0 \leq k \leq n : k \in B\}$$ exists, it is called
the asymptotic frequency of $B$. A partition of $\N$ can be
thought of as a sequence $B_1,B_2,...$ of disjoint blocks whose
union is $\N$. The labeling obey the following rule : if $B_i$ is
not empty, then its least element is $i$. Call $\P$ the space of
the partitions of $\N$, and recall that $\P$ is a metric compact
space, see \cite{kingman}.

A finite permutation $\sigma$ (i.e. a bijection $\N \rightarrow
\N$ such that $\sigma(n)=n$ for $n$ large enough) acts on a
partition $\pi$ in the following way : for any $i$ and $j$ in
$\N$, $i$ and $j$ are in the same block of $\sigma(\pi)$ iff
$\sigma(i)$ and $\sigma(j)$ are in the same block of $\pi$, this
equivalence relation can be identified as a partition and thus
completely define $\sigma(\pi)$.

A measure $\mu$ on $\P$ is said \textit{exchangeable} if for any
measurable set $A \subseteq \P$, for any finite permutation
$\sigma$
$$\mu(A)=\mu(\sigma(A)),$$ where $\sigma$ acts on the sets in the
obvious way.

A $\P$-valued process $\Pi$ is said exchangeable if the permuted
process $\sigma(\Pi)$ has the same distribution as the original
process $\Pi$. For instance the $\P$-valued process $\Pi_F(t)$
presented in the introduction is exchangeable.

For all $B \subseteq \N$, let $P_B$ be a probability on the
partitions of $B$. For all $\pi=(B_1,B_2,..) \in \P$, let
$P_{\pi}$ be the distribution of the partition with blocks
$B_{(1,1)},B_{(1,2)},.....,B_{(2,1)},B_{(2,2)},...$ where
$\pi^{(i)}=(B_{(i,1)},B_{(i,2)},...)$ is a partition of $B_i$ and
has law $P_{B_i}$. The family $(P_{\pi} , \pi \in \P)$ is, in the
terminology of Pitman \cite{pitman}, a fragmentation kernel on
$\P$.

\begin{definition}
Call $\P$-fragmentation any exchangeable $\P$-valued Markov
process, starting from the trivial partition ($\N$ is the only
non empty block), which is continuous in probability and has
fragmentation kernels as its transition semi-group.
\end{definition}

We briefly recall some definitions and results on
$\P$-fragmentations. If $\pi$ is a random exchangeable partition,
by a result of Kingman \cite{kingman} (see also Aldous
\cite{aldous} for a simpler proof),  every block of $\pi$ has an
asymptotic frequency almost surely, i.e. $|B_i|$ exists with
probability 1 for all $i=1,...$.

We call an exchangeable $\P$-valued process $\Pi$ \textit{nice} if
with probability 1, $\Pi(t)$ has asymptotic frequencies for all $t
\geq 0$ simultaneously. Evans and Pitman \cite{evans-et-pitman}
have shown that it is always the case when $\Pi$ is an
exchangeable $\P$-process with proper frequencies (i.e. for each
$t\geq 0$, $\sum_{i \in\N}|B_i(t)|=1$ almost surely ), and Bertoin
\cite{bertoin1} proved that so-called homogeneous fragmentation
were nice. Observe that when $\Pi(t)$ is nice, the ordered
sequence of the asymptotic frequencies is well defined and is a
$\sfleche$-valued process.

As we shall construct a Markovian semi-group on $\sfleche$, we
need a notion slightly more general than the asymptotic frequency,
well defined for any subset $B$ of $\N$. We write
$$\Lambda(\Pi(t))=
(\Lambda_1(\Pi(t)),\Lambda_2(\Pi(t)),...)^{\downarrow}=(\lambda_1(t),\lambda_2(t),..)$$
for the decreasing rearrangement of the quantities
$$\Lambda_i(\Pi(t))=\liminf_{n \rightarrow \infty} \frac{1}{n}
\# \{ k\leq n : k \in B_i(t)\} .$$ By extension we also note
$$\Lambda(B)=\liminf_{n \rightarrow \infty} \frac{1}{n} \#
\{ k\leq n : k \in B\}$$ for any $B \subseteq \N$.

$\Lambda$ is a functional of $\Pi(t)$ that takes its values in
$\sfleche$. We stress that $\Lambda$ is not continuous.

Next for every $C \subseteq \N$ and every $\pi = \{B_1,B_2,...\}
\in \P$, we define the partition of $C$ induced\footnote{there is
in fact another natural way of defining this partition : it is to
take the image of $\pi$ by the mapping that sends $\N$ onto
$C=\{c_1,c_2,...\}$ (where $c_1<c_2<...$) :
$$\pi \circ C =\left( \{ c_j : j\in B_i \}_{ i=1,...}\right)$$ Suppose now that $\pi$ is an exchangeable random $\P$-valued
variable, for all $k>0$, for any finite permutation $\sigma$ such
that $$\forall i \leq k  ; \sigma(i)=c_i ,$$ $\pi$ and
$\sigma(\pi)$ have same law, thus in the sense of the equality of
the finite-dimensional margins $\pi \circ C$ and $\pi \cap C $
have same law. Thus in fact any definition could be taken
indifferently.} by $\pi$ :
$$\pi \cap C =\left( B_1 \cap C , B_2 \cap C,... \right).$$

\begin{definition}
A $\P$-fragmentation $\Pi=(\Pi(t),t\geq 0)$ is called self-similar
with index $\alpha \in \R$ if :
\begin{enumerate}
\item $\Pi$ starts a.s. from the trivial partition.
\item The ranked fragmentation $\Lambda(\Pi)$ associated to $\Pi$
is continuous in probability.
\item For every $B \subseteq \N$, $\forall
t \geq 0$ $P_B(t)$ (in the above notations) is the distribution of
$\Pi(t\Lambda(B)^{\alpha}) \cap B$.
\end{enumerate}
When $\alpha=0$ we will say that $\Pi$ is a homogeneous
fragmentation.
\end{definition}

Following Kingman \cite{kingman} (see also \cite{aldous} for a
survey), to each $s=(s_1,s_2,...) \in \sfleche$ one can associate
a unique exchangeable probability measure $\mu_s$ on $\P$ such
that $\mu_s$-almost every partition has ranked asymptotic
frequencies $s$.

This is how one proceeds : let $(X_i)_{i\in \N}$ a family of iid
variables such that  $\forall k \in \N , P(X_i=k)=s_k$ and
$P(X_i=-i)=1-\sum_k s_k$, then define the $s$ -paintbox
\footnote{The reason for the name (due to Kingman) is the
following : imagine that we have a choice of colors $(c_k)_{k \in
\N}$. Then paint each integer $n$ independently with a randomly
chosen color, $c_k$ with probability $s_k$. Then the partition of
$\N$ defined by the equivalence relation "being of the same color"
is the $s$-paintbox process.} partition ( or "$s$-paintbox
process") $\Pi$ by the equivalence relation
$$\forall i,j \in \N , i \sim j \Leftrightarrow X_i=X_j.$$ We
denote by $\mu_s$ the law of the $s$-paintbox process. It is clear
by the LLN that $\mu_s$-almost surely $\Lambda(\Pi)=s$.

For each self-similar $\P$-fragmentation one can take the
associated $\sfleche$ ranked fragmentation, thus defining a map
from $\P$-fragmentation laws into $\sfleche$-fragmentation laws.
Suppose now that $\Pi_1$ and $\Pi_2$ are two self-similar
$\P$-fragmentations such that for any fixed $t$ the $\sfleche$
variables $\Lambda(\Pi_1(t))$ and $\Lambda(\Pi_2(t))$ have same
law. $\Pi_1(t)$ and $\Pi_2(t)$ being exchangeable, by de Finetti's
theorem (see \cite{aldous} one can show that they are mixture of
paintbox processes directed respectively by $\Lambda(\Pi_1)$ and
$\Lambda(\Pi_2)$, i.e. $$P(\Pi_{1,2} \in A)=\int_{\sfleche}
\mu_s(A) P(\Lambda(\Pi_{1,2}) \in ds).$$ We conclude that they
have the same distribution. So to every $\P$-fragmentation
corresponds a different $\sfleche$-fragmentation. Our first result
will be to show that there is in fact a one to one relation.

\subsection{From Ranked to Partition Fragmentations}

Let $\Pi$ be a nice self-similar fragmentation of index $\alpha$,
then it is not difficult to show that its asymptotic frequencies
$\Lambda(\Pi)$ form a self-similar ranked fragmentation of index
$\alpha$. Conversely we shall now show that to each
$\sfleche$-fragmentation $\lambda$ we can associate a
$\P$-fragmentation $\Pi$ with same index of self-similarity such
that $\lambda=\Lambda(\Pi)$.

\medskip
\begin{proposition} \label{relations ranked partition}
We have the following relations between  $\sfleche$ and $\P$
fragmentations :
\begin{enumerate}
\item If $\Pi$ is a $\P$-fragmentation then $\Lambda(\Pi)$
has the finite-dimensional marginal distributions of an
$\sfleche$-fragmentation. Moreover $\Lambda$ preserves
self-similarity.
\item if $\lambda$ is a $\sfleche$-fragmentation, then we
can construct $\Pi_{\lambda}$ an exchangeable $\P$-fragmentation
such that $\Lambda(\Pi_{\lambda}) \eqlaw \lambda$. Moreover this
construction preserves self-similarity.

\end{enumerate}

\end{proposition}
\medskip

The first point is clear, the difficulty here lies in the second
part of this proposition. The main idea is that as $\P$ is a
compact metric space, it is enough to construct an adequat
Markovian semi-group to ensure the existence of the desired
$\P$-process. Then the conservation of the index will be a simple
consequence of our construction.

Let $(P_t(S),t \geq 0, S \in \sfleche)$ be a transition kernel on
$\sfleche$ generated, in the notation of definition (\ref{def1}),
by the family $(P_t(l),t \geq 0, l \in [0,1])$.

Let $\tilde{P}_t(l)$ be the image of $P_t(l)$ by $g_{l^{-1}}$ the
map $ (x_1,x_2,...) \rightarrow (x_1/l,x_2/l,...)$. Let
$(Q_t(l,.),l \in [0,1], t \geq 0)$ be a family of probability
measures on $\P$ where, for a fixed $t$, $Q_t(l)$ is a mixture of
$s$-paintbox processes directed by $\tilde{P}_t(l)$, i.e. for $A
\subseteq \P$
$$Q_t(l,A)=\int_{\sfleche} \mu_s(A) \tilde{P}_t(l,ds).$$

For $B\subseteq \N$ define $Q_t(B)$ the distribution of $\Pi_B
\cap B$ where $\Pi_B$ is a $\P$ valued random variable with law
$Q_t(\Lambda(B))$. Practically this means that one begins by
drawing a variable $\lambda_B$ with law $\tilde{P}_t(\Lambda(B))$
and then the $\lambda_B$-paintbox partition and then take its
intersection with $B$.

Now let $\pi=(\pi_1,\pi_2,...) \in \P$ and $\forall t \geq 0$ let
$(\Pi_{\pi_i}(t))_{i \in \N}$ be a sequence of independent
variables with respective law $Q_{t}(\pi_i)$. Define $Q_{t}(\pi)$
the law of the partition whose blocks are the blocks of the
$(\Pi_{\pi_i}(t),i \in \N)$.

Our proof of Proposition \ref{relations ranked partition} shall
thus consist in showing that the family $(Q_{t}(\pi),\pi \in \P,
t\geq 0)$ forms a semi-group.

\begin{proof}
From the above description it should be clear that it suffices to
show
\begin{eqnarray}\label{semi-group}
\forall \pi \in \P, Q_{t+u}(\pi)=\int_{\pi' \in \P}
Q_{t}(\pi')Q_{u}(\pi,d\pi')
\end{eqnarray}
in the obvious notation. If for any subset $B$ of $\N$ we note
$\P_B$ for the space of the partitions of $B$, by construction,
(\ref{semi-group}) is equivalent to
\begin{eqnarray}\label{semi-groupbis}
\forall B \subseteq \N, Q_{t+u}(B)=\int_{\pi' \in \P_B}
Q_{t}(\pi')Q_{u}(B,d\pi').
\end{eqnarray}

We can reformulate (\ref{semi-groupbis}) as : $Q_{t+u}(B)$ is the
distribution of the random partition $\Pi(t,u)$ of $B$ (and this
is what we actually shall prove) obtained by the following
two-steps procedure :
\begin{enumerate}
\item draw $\Pi(u)=(\pi_1(u),\pi_2(u),...)$ an exchangeable partition of $B$ with
law $Q_{u}(B)$.
\item given $\Pi(u)$ draw a sequence $\left(\Pi_{\pi_i(u)}(t)\right)_{i\in
\N}$ of independent $\P_{\pi_i(u)}$-variables with respective law
$Q_t(\pi_i(u))_{i \in \N}$.
\item $\Pi(t,u)$ is just the collections of all the blocks of the
$\Pi_{\pi_i(u)}(t)$.
\end{enumerate}

We begin by proving so for $B=\N$. By construction we can always
suppose that $\Pi(u)$ is a mixture of paintbox processes directed
by $P_u((1,0,...))$, i.e. conditionally on $\lambda(u)$ (a random
variable with law $P_u((1,0,...))$), $\Pi(u)$ is a
$\lambda(u)$-paintbox process (resp. for each $i \in \N$
$\Pi_{\pi_i(u)}(t)$ is constructed by taking the intersection of
$\pi_i(u)$ and a $\lambda^{(i)}$ paintbox-process where
$\lambda^{(i)}$ is a $\sfleche$-variable with law
$\tilde{P}_t(|\pi_i(u)|)$.)

This means that conditionally on
$\lambda(u)=(\lambda_1(u),\lambda_2(u),...)$ one draws an i.i.d.
sequence of variables $(X_i)_{i \in \N}$ with values in $\N$ whose
law is $P(X_1=k)=\lambda_k(u)$ for any $k \geq 1$ and
$P(X_i=-i)=1-\sum_n \lambda_n(u)$ which determines $\Pi(u)$ (idem
with each $\Pi_{\pi_i(u)}(t)$ denoting $(Y^{(i)}_k)_{k \in \N}$
the appropriate sequence of variables).

Fix $\phi$ a bijection from $\N^2$ in $\N$ and define the
coordinate $\phi^{-1}(k)=(\alpha_k,\beta_k)$ for all $k$, then
for $i$ and $j$
\begin{eqnarray*}i \overset{\Pi(u,t)}{\sim} j
& \Leftrightarrow &  \{X_i=X_j \text{ and }
Y^{(X_i)}_i=Y^{(X_j)}_j\}
\\ & \Leftrightarrow & \{(X_i,Y^{(X_i)}_i)=(X_j,Y^{(X_j)}_j)\}
\\ & \Leftrightarrow & \phi(X_i,Y^{(X_i)}_i)=\phi(X_j,Y^{(X_j)}_j)
\\ & \Leftrightarrow & Z_i=Z_j
\end{eqnarray*}
where $Z$ is an obvious notation $Z_i=\phi(X_i,Y_i^{(X_i)})$.

Then we have the identity
$$\{Z_i=k\}=\{X_i=\alpha_k,Y^{(\alpha_k)}_i=\beta_k\}$$ As
conditionally on $(\lambda^{(i)}(t))_{i\in \N}$ and $\lambda(u)$
the $X_i$ are i.i.d. as well as the sequences $(Y^{(i)}_k)_{k \in
\N}$ and are all independent between them, we see that the $Z_i$
are also i.i.d. As $$ i \overset{\Pi(t,u)}{\sim}j \Leftrightarrow
\{ Z_i=Z_j \}$$ $\Pi(t,u)$ is exchangeable.

The law of an exchangeable random partition is completely
determined by the law of its asymptotic frequencies, here the
$\lambda_i(u)\times \lambda^{(i)}(t).$ As $\lambda(.)$ is a
$\sfleche$-fragmentation we have by construction that
$$ \left( (\lambda_i(u)\times \lambda^{(i)}(t))_{i \in \N}
\right)^{\downarrow} \eqlaw \lambda(t+u).$$  So $\Pi(t,u)$ has law
$Q_{t+u}(\{\N\})$.

Then take $B$ a subset of $\N$. By construction $Q_{t+u}(B)$ is
the law of $\tilde{\Pi}_{t+u}(B) \cap B$ where
$\tilde{\Pi}_{t+u}(B)$ is a certain $\P$-variable and
$\Pi(t,u)=\tilde{\Pi}(t,u) \cap B$ where $\tilde{\Pi}(t,u)$ is a
certain variable. It is clear that replacing the generating
family $(Q_t(l), t \geq 0 , l \in [0,1])$ by
$(Q'_t(l)=Q_t(\Lambda(B)l), t \geq 0 , l \in [0,1])$ the above
arguments yields $\tilde{\Pi}(t,u) \eqlaw \tilde{\Pi}_{t+u}(B)$
and thus for all $B \subseteq \N$ $$\Pi(t,u) \eqlaw
\Pi_{t+u}(B).$$

So we have proved the existence of a Markov $\P$-process $\Pi$
with semi-group $Q_{\pi}(t)$, which, by construction, is a
fragmentation whose asymptotic frequencies has same distribution
(in the sense of finite-dimensional distributions) as $\lambda_t$
our starting $\sfleche$-fragmentation.

For each  ranked fragmentation $\lambda$ we can thus construct a
partition fragmentation $\Pi_{\lambda}$ such that
$\Lambda(\Pi_{\lambda})$ has same law as $\lambda$.

We now turn to the conservation of self-similarity : suppose
$\lambda$ is a self-similar $\sfleche$-fragmentation with index
$\alpha$, so $\tilde{P}_l(t)=P_1(tl^{\alpha})$, looking at the
above construction of the semi-group of $\Pi_{\pi}$ shows that
$Q_t(l)=Q_{1}(tl^{\alpha})$, so $\Pi$ is also self-similar of
index $\alpha$.
\end{proof}

It is now natural to look for some explicit construction of ranked
fragmentation i.e. an equivalent of Theorem 1 in \cite{bertoin1}.

\section{Homogeneous fragmentation}

In \cite{bertoin1} J. Bertoin shows how a homogeneous
$\P$-fragmentation process can be decomposed into a Poisson point
process of partitions, whose distribution is determined by the
so-called  characteristic measure. We will begin by recalling the
facts we need on this topic, and then tackle the analog problem
for $\sfleche$ homogeneous fragmentations.

\subsection{L\'evy-It\^o decomposition of homogeneous \\ $\P$-fragmentations}

The distribution of a homogeneous $\P$-fragmentation $\Pi$ is
determined by an exchangeable measure $\kappa$ on $\P$, called the
characteristic measure of $\Pi$, that assign zero mass to the
trivial partition and verifies the condition
$\kappa(\P^*_2)<\infty $ where $\P^*_2$ is the set of the
partitions of $\N$ for which $1$ and $2$ does not belong to the
same block. Given such a measure $\kappa$, one can construct an
homogeneous $\P$-fragmentation admitting $\kappa$ as its
characteristic measure as follows : Let $K=((\Delta(t),k(t)),t\geq
0)$ a Poisson point process with values in $\P \times  \N$ with
intensity measure $M:=\kappa \otimes \#$ where $\#$ stands for the
counting measure on $\N$. This means that for a measurable set
$A\subseteq \P \times \N$ with $M(A) < \infty$, the counting
process $$N^A(t)=\text{Card}(s \in [0,t]: (\Delta(s),k(s)) \in A),
t\geq 0)$$ is a Poisson process with intensity $M(A)$, and to
disjoint sets correspond independent processes.

Then one can construct a unique $\P$-valued process
$\Pi_{\kappa}=(\Pi_{\kappa}(t),t\geq 0)$ started from the trivial
partition, with c\`adl\`ag sample paths, such that $\Pi_{\kappa}$
only jumps at time $t$ at which $K$ has an atom
$(\Delta(t),k(t))$, and in that case $\Pi_{\kappa}(t)$ is the
partition whose blocks are the $B_i(t_-)$ (the blocks of
$\Pi_{\kappa}(t_-)$) except for $B_{k(t)}(t_-)$ which is replaced
by the partition of $B_{k(t)}(t_-)$ induced by $\Delta(t)$ (that
is $\Delta(t) \cap B_{k(t)}(t_-)$).

$\Pi_{\kappa}$  is a homogeneous $\P$-fragmentation with
characteristic measure $\kappa$. Conversely, any homogeneous
$\P$-valued fragmentation $\Pi$ has the same law as $\Pi_{\kappa}$
for some unique exchangeable measure $\kappa$.

As a consequence of Kingman's representation of exchangeable
partitions \cite{kingman}, every exchangeable partition measure
can be decomposed as the sum of a \textit{dislocation} measure and
an \textit{erosion} measure :
\begin{itemize}
\item  $\delta_{\pi}$ stands for the Dirac point mass at $\pi \in \P$, for all $n \in
\N$ let $\epsilon_n$ be the partition of $\N$ with only two
non-voids blocks : $\{n\}$ and $\N \backslash \{n\}$, then for
every $c \geq 0$, the measure
$$\mu_c=c \sum_{n=1}^{\infty}\delta_{\epsilon_n}$$ is an
exchangeable measure. The $\mu_c$'s are called \textit{erosion}
measures.
\item The dislocation measures are constructed from so-called L\'evy
measures on $\sfleche$. We call a measure $\nu$ on $\sfleche$ a
L\'evy measure if $\nu$ has no atom at $(1,0,0,..)$ and verifies
the integral condition
$$\int_{\sfleche} (1-s_1)\nu(ds) < \infty$$ where
$s=(s_1,s_2,...)$ denotes a generic sequence in $\sfleche$. The
mixture of paintbox processes
$$\mu_{\nu}(.)=\int_{\sfleche} \mu_s(.) \nu(ds)$$ is a
measure on $\P$, called the \textit{dislocation} measure directed
by $\nu$.
\end{itemize}

Then for any $\kappa$ exchangeable partition measure there exists
a unique $c \geq 0$ and a unique L\'evy measure $\nu$ such that
$\kappa= \mu_c + \mu_{\nu}$.

Thus the law of a homogeneous $\P$-fragmentation is completely
characterized by the pair $(\nu,c)$. Using Proposition
\ref{relations ranked partition}, we conclude that :

\begin{corollary}
There is a bijective correspondence between the laws of
homogeneous ranked fragmentations and the pairs $(\nu,c)$ where
$\nu$ is a L\'evy measure on $\sfleche$ and $c \geq 0$.
\end{corollary}

A ranked fragmentation is thus completely characterized (in terms
of distribution) by the pair $(\nu,c)$ associated to its law.

We would like to transfer the Poisson point process construction
of $\P$-fragmentations to $\sfleche$-fragmentations. The main
difficulty in doing so comes from the lack of a genealogy
structure in this new setting.

To illustrate this, let $K=(\Delta(t),k(t))$ a PPP with measure
intensity $\mu_{\nu} \times \#$ and $\Pi$ the corresponding $\P$
fragmentation (hence with no erosion), and suppose that at time
$t$ the $k$-th block of $\Pi(t_-)$ (i.e. its least element is $k$)
fragments, or otherwise said
$\Lambda(B_k(t_-))>\Lambda(B_k(t))>0$. Then it is clear that at
time $t$ there is also a dislocation in the associated ranked
fragmentation $\lambda=\Lambda(\Pi)$. The label of the mass of
$\lambda(t_-)$ that fragments, noted $\Phi(t_-,k)$, is an integer
that depends on $\Pi(t_-)$ and $k$ and can informally be seen as
the rank of the size of the $k(t)$-th block of $\Pi(t_-)$. In the
same way that $\Pi$ is constructed from $K$, one might hope that
$\Lambda(\Pi)$ is constructed from $(\Lambda(\Delta_t),
\Phi(t_-,k_t))$ but we will still have to show that this last
point process is a Poisson point process with the right intensity,
then that the jump-times of $\lambda$ are exactly the atom times
of $K$ and finally that $\lambda$ is a pure-jump process (in a
sense to be defined).

But first we show how to get rid of erosion.

\subsection{Erosion in homogeneous ranked fragmentation}

Let us first examine the trivial case when the  fragmentation is
pure erosion. It is then intuitively clear that the homogeneity in
time and space entails that the ranked fragmentation $\lambda(t)$
with values in $S^{\downarrow}$ with characteristics $(0,c)$
(where the 0 means that the measure $\nu$ is trivial with mass 0)
is given by $$\lambda(t)=(e^{-ct},0,0,...)$$ To demonstrate this
define $$k=\mu_c$$ with $c>0$, and let $\Pi$ be the
$\P$-fragmentation associated to the P.P.P. $K=(\Delta(t),k(t))_{t
\geq 0}$ with intensity $\mu_c \otimes \# $ and values in
$\P\times \N$. $\Pi$  can be thought of as an isolation process,
indeed at each jump time of $K$, say $t$, some point of $\N$, say
$n$, is designated, (i.e. $\Delta(t)=\delta_{\epsilon_n}$). If the
block containing $n$, $\beta(n,t)$, is not reduced to the
singleton $\{n\}$, then it is fragmented into $\{n\}$ and
$\beta(n,t) \backslash \{ n \} $, "$n$ is isolated from its
block", else nothing happens. Hence, at all time there is only one
block which is not a singleton, by an argument that will be
established thereafter in Theorem \ref{prop saut pur}, we can
always suppose that this block also contains $1$. If we consider
the restriction of $\Pi$ to $\{1,2,..,n\}$, denoted by
$\Pi^{(n)}$, then $\Pi^{(n)}$ only jumps at atom-times of $K$ for
which $k_t=1$ and $\Delta_t \in
\{\delta_{\epsilon_1},\delta_{\epsilon_2},...,\delta_{\epsilon_n}\}$.
The restriction of the Poisson Process to this set is a Poisson
Process with intensity of finite mass and have thus discrete
jump-times. The processes of the times of exclusion of each point
are independent one of the other. By standard calculation on
Poisson processes the probability that a given point have been
excluded at time $t$ is $\exp(-tc)$, thus the law of the number of
point excluded at time $t$ is a Bernoulli with parameter
$(e^{-ct})$. By the law of the large number, at time $t$, the
asymptotic frequency of the only block not reduced to a singleton
is $(e^{-ct})$ almost surely. So a.s. for every $t \in \mathbb{Q}$
$$\Lambda(\Pi_t)=((e^{-ct}),0,0,...)$$ and as $\lambda_1(t)$ is
monotone decreasing the relation holds almost surely for all $t$.
This result is the key for the following.

\begin{proposition}
If $\tilde{\lambda}$ is a homogeneous $(\nu,0)$ ranked
fragmentation, then $\lambda=(e^{-ct} \tilde{\lambda}(t), t \geq
0)$ is a homogeneous $(\nu,c)$ ranked fragmentation.
\end{proposition}
\begin{proof}
Let $\tilde{\Pi}$ and $\Pi$ be some homogeneous partition
fragmentations with characteristics $(\nu,0)$ and $(\nu,c)$
respectively. Then call $\tilde{\lambda}$ the process of the
ordered asymptotic frequencies of $\tilde{\Pi}$ and $\lambda$
those of $\Pi$. Suppose $\Pi$ is constructed on the Poisson point
process $K=(\Delta(t),k(t) , t \geq 0)$ with characteristic
measure $\mu_{\nu} + \mu_c$. Let $K_1=(\Delta(t),k(t) , t \geq
0)$ the Poisson point process with characteristic measure
$\mu_{\nu} \otimes \# $ and $K_2=(\Delta(t) , t \geq 0)$ the
Poisson point process with characteristic measure $\mu_{c} $.

Thus $\Pi$ appears as (i.e. is equal in law to) the intersection
of $\Pi_1$ (constructed on $K_1$) and a pure erosion process
$\Pi_2$ (constructed on $K_2$), i.e. $\Pi= \Pi_1(.) \cap
\Pi_2(.)$ defined by the equivalence relation
$$\forall i,j \in \N : (i \overset{\Pi_1(.) \cap \Pi_2(.)}{\sim}
j) \Leftrightarrow \left((i \overset{ \Pi_1(.)}{\sim} j)\text{ and
} (i \overset{ \Pi_2(.)}{\sim} j)\right).$$

Given a random exchangeable subset of $\mathbf{N}$, say $A$,
independent of $(\Pi_2(t))_{t \geq 0}$, with random asymptotic
frequency $l$, the asymptotic frequency of the subset of $A$
defined as the points that have not been excluded up to time $t$
is $le^{-ct}$ a.s. for all $t$.

Therefore $\Lambda(\Pi(t)) \eqlaw e^{-ct}\Lambda(\Pi_1(t))$. As we
can always suppose that a ranked fragmentation is the associated
ranked fragmentation of some partition fragmentation the result is
proven.
\end{proof}

Thus it suffices to know how to construct a homogeneous ranked
fragmentation without erosion from a PPP to know how to construct
any homogeneous ranked fragmentation.

\subsection{Construction of homogeneous ranked fragmentation with no erosion}

Let $\lambda$ be an $\sfleche$-fragmentation, with characteristics
$(\nu,c)$, then for every $k \in \N$ the process
$\lambda_1(t)+...+\lambda_k(t)$ is monotone decreasing. $\lambda$
is said to be a pure jump process if for any $k$,
$\lambda_1(t)+...+\lambda_k(t)$ is a pure jump process.

In the following we shall focus on the case where for each fixed
$t$ there is a infinite number of fragments of strictly positive
size almost surely. A necessary and sufficient condition for this
is $$\nu(s \in \sfleche : s_2>0)=\infty.$$ Indeed, fix $t>0$ and
suppose that $\lambda_1(t)>0$. Then, for any $\epsilon>0$, during
the time interval $[t-\epsilon,t]$, $\lambda_1$ has been affected
by an infinite number of dislocation such that at least one small
fragment detached from the main one, thus an infinite number of
fragments have been created, and the life-time of those variables
form a sequence of independent identically distributed random
variables, thus with probability one an infinite number of them
have survived at time $t$. The same line of arguments also shows
that $\inf \{ t \geq 0 : \lambda_1(t)=0\}=\infty$ almost surely.

Although most of the following results are still true for any
$\nu$, making this hypothesis enables us to focus on the most
interesting case and to avoid some technical difficulties.

\begin{theorem}\label{prop saut pur}
Let $\lambda$ be a homogeneous $\sfleche$-fragmentation with no
erosion ($c=0$) and L\'evy measure $\nu$ as above (i.e. $\nu(\{s
:s_2>0 \}=\infty)$. Then
\begin{enumerate}
\item $\lambda$ is a pure jump process.
\item there exists a PPP $K=(S(t),k(t))_{t\geq 0}$ with
values in $\sfleche \times \N$ and intensity measure $\nu \otimes
\#$, such that the jumps of $\lambda$ correspond to the atoms of
$K$. More precisely, $\lambda$ only jumps at times at which
$(S(t),k(t))$ has an atom, and at such a time $\lambda(t)$ is
obtained from $\lambda(t_-)$ by dislocating the $k(t)$-th
component of $\lambda(t_-)$ by $S(t)$ (i.e. replacing
$\lambda_{k(t)}(t_-)$ by the sequence $\lambda_{k(t)}(t_-) S(t)$)
and reordering the new sequence of fragments. Conversely if
$(S(t),k(t))$ is an atom then $\lambda$ has a jump at $t$, i.e.
$\lambda_i$ jumps at $t$ for some $i$.
\end{enumerate}
\end{theorem}

Although this result is intuitive in regard to the equivalence
between $\P$ and $\sfleche$ fragmentation, it requires some
technical work.

We give ourselves a homogeneous $(\nu,0)$ $\sfleche$-fragmentation
$\lambda$ with $\nu$ verifying $\nu(s_2>0)=\infty$. There is no
loss of generality in  supposing that $\lambda$ is constructed as
follows : Call $H=((\Delta(t),k(t)))_{t \geq 0}$ a PPP with
measure intensity $\mu_{\nu} \otimes \#$ with values in $\P
\times \N$. Let $\Pi$ be the homogeneous $(\nu,0)$
$\P$-fragmentation constructed on $H$, then define
$$\lambda=\Lambda(\Pi).$$ Call $\ft=\sigma\{\Pi_s , s\leq
t\}$ the natural filtration of the $\P$-fragmentation $\Pi$.

Then at any time $t$, call $\phi(t,.)=\phi_t(.)$ the random,
$\ft$ measurable application from $\N \rightarrow \N \cup \infty$
(where $\infty$ serves as a cemetery point) defined as
\begin{itemize}
\item if $|B_k(t)|>0$ then  $\phi(t,k)$ is the rank of the asymptotic frequency of
$B_k(t)$ (it is well defined because the number of blocks of
greater asymptotic frequencies is always finite with an upper
bound of $|B_k(t)|^{-1},$ and in case two blocks have the same
asymptotic frequency, they are ranked as their least element).
\item if $|B_k(t)|=0 $ (with the convention $|\o|=0$) then $\phi(t,k)=\infty$
\end{itemize}
We also note $\tilde{k}(t)=\phi(t_-,k(t))$. Note that under our
hypothesis that there is always an infinite number of fragments
$\forall t \geq 0, \N \subset \{ \phi(t,k) , k \in \N \}$.

We will first prove that the point process image of $H$, noted
$\tilde{K}$, whose atoms are the points of
$(\Lambda(\Delta(t)),\tilde{k}(t))_{t \geq 0}$ such that
$\tilde{k}(t)\in \N$, is a Poisson point process with measure
intensity $\nu \otimes \#$. Then we will show that this is also
the process of the jumps of $\Lambda(\Pi)$ and this last process
is a pure jump process so it can wholly be recovered from
$(\Lambda(\Delta(t)),\tilde{k}(t))_{t \geq 0}$. This will complete
the proof of Theorem \ref{prop saut pur}.

\begin{lemma}\label{poisson}
The point process $\tilde{K}(t)$ derived from
$(\Lambda(\Delta(t)),\tilde{k}(t))_{t \geq 0}$ by only keeping the
atoms such that $\tilde{k}(t) \neq \infty$ is a Poisson point
process with intensity measure $\nu \otimes \#$;
\end{lemma}

\begin{proof}
Let $A$ be a subset of $\sfleche$ such that $\nu(S)<\infty$. For
$i=1,...$ let $$N_A^{(i)}(t)=\#\{u \leq t : \Lambda(\Delta(u)) \in
A , k(u) = i\} $$  Then set $$N_A(t)=\#\{u \leq t :
\Lambda(\Delta(u)) \in A , \tilde{k}(u) = 1\}. $$  $N_A(t)$ is
increasing, right-continuous with left-limits with jumps of size 1
(the $N_A^{(i)}(t)$ being independent Poisson processes they do
not jump at the same time almost surely).  By definition we have
\begin{eqnarray*}
dN_A(t) &=& \sum_{i=1}^{\infty}\mathbf{1}_{\{\phi(t_-,i)=1\}}
dN_A^{(i)}(t)
\end{eqnarray*}

Define
\begin{eqnarray*}
d\tilde{N}_A^{(i)}(t)=\mathbf{1}_{\{\phi(t_-,i)=1\}}dN_A^{(i)}(t).
\end{eqnarray*}
It is clear that $\mathbf{1}_{\{\phi(t_-,i)=1\}}$ is adapted and
left-continuous in $(\ft)$ and hence predictable. The
$N_A^{(i)}(.)$ are i.i.d. Poisson processes with intensity
$\nu(A)$ in $(\ft)$. Thus, for each $i$ the process
$$M_A^{(i)}(t)=\tilde{N}_A^{(i)}(t)-\nu(A)\int_0^t
\mathbf{1}_{\{\phi(u_-,i)=1\}}du=\int_0^t\mathbf{1}_{\{\phi(u_-,i)=1\}}d(N_A^{(i)}(u)-\nu(A)u)$$
is a square integrable martingale. Then define
$$M_A(t)=\sum_{i=1}^{\infty}
\int_0^t\mathbf{1}_{\{\phi(u_-,i)=1\}}d(N_A^{(i)}(u)-\nu(A)u)$$

Note $f_i(t)=\mathbf{1}_{\{\phi(t_-,i)=1\}}$, then, for all $i
\neq j , \forall t \geq 0 \text{ , } f_i(t)f_j(t) = 0$, and $
\forall t \text{ , } \sum_{i=1}^{\infty}f_i(t)=1$.

As the $N_A^{(i)}(t)$ are independent Poisson processes they do
not jump simultaneously and so the martingales $M_A^{(i)}(t)$ do
not either. They are thus orthogonal (see for example chapter 8,
Theorem (43)-D in \cite{dellacherie-meyer} for a proof). Moreover
the oblique bracket of $M$ is

\begin{eqnarray*}
\langle M_A \rangle (t)
&=&\sum_{i=1}^{\infty}<\int_0^tf_i(u)d(N_A^{(i)}(u)-\nu(A)u)>
\\ &=& \nu(A)t
\end{eqnarray*}
So $M_A$ is a $L_2$ martingale.

So we have demonstrated that $N_A(t)$ is increasing,
right-continuous, left limited with jump of size 1 with
compensator $\nu(A)t$. Using classical results (see  for instance
chapter 2.6 in \cite{ikeda-watanabe}, Theorem 6.2) we conclude
that $N_A(t)$ is a Poisson process with intensity $\nu(A)$. Now
take $B \in \sfleche$ such that $A \cap B =\O$, we can use the
same construction as above replacing $A$ with $B$ and the fact
that $N_A^{(i)}(t)$ and $N_B^{(i)}(t)$ are independent Poisson
processes in the same filtration to see that
$$\{\Lambda(\Delta(u))  : u \geq 0 , \tilde{k}(u)=1\} $$
is a P.P.P. with intensity measure $\nu$. the same arguments yield
that
$$(\{ \Lambda(\Delta(u)) : u \geq 0 , \tilde{k}(u)=2\} $$
is also a P.P.P. with measure intensity $\nu$. It is clear that
$N_1$ and $N_2$ have no jumps in common because the
$N_A^{(i)}(t)$'s does not, so they are independent. By iteration
this show that $(\Lambda(\Delta(t)),\tilde{k}(t))$ is a P.P.P.
with measure intensity $\nu \times \#$.
\end{proof}

Let $K$ be a P.P.P. on $\P \times \N$ with intensity measure
$\mu_{\nu} \otimes \#$ and
$\Pi=(\Pi(t),t\geq0)=((B_1(t),B_2(t),...),t\geq 0)$ the $(\nu,0)$
$\P$-fragmentation constructed from $K$, and define
$\lambda=\Lambda(\Pi)=(\lambda_1(t),\lambda_2(t),..)$ the ordered
vector of asymptotic frequencies. In the case considered here $\Pi
$ is nice so  almost surely for all $t$ $|B_i(t)|$ exists for all
$i \in \N$. Recall that $\phi(t,k)$ is the rank of the asymptotic
frequency $|B_k(t)|$ at time $t$.

We now need to show that $\lambda$ is a pure jump process in the
sense that for each $k$ the decreasing process
$\lambda_1+...+\lambda_k$ is pure jump and that all his jumps are
indeed images of some atoms of $K$ ($\Lambda$ being not continuous
it is not \textit{a priori} evident).

In \cite{bertoin1} it is shown that $|B_1(t)|$, the asymptotic
frequency of the block that contains $\{1\}$, is the inverse of
the exponential of a subordinator with $0$-drift, and so it is a
pure-jump process. By the Markov and homogeneity property this
implies that for all $i>1$ the process $|B_i(t)|$, the asymptotic
frequency of the block that contains $i$, is c\`adl\`ag, started
at $0$, such that at $\tau_i=\sup\{t \geq 0 : |B_i(t)|=0 \}$ we
have $|B_i(\tau_i)|>0$ (i.e. it leaves 0 with a jump), and after
$\tau_i$ the process $\frac{|B_i(t-\tau_i)|}{|B_i(\tau_i)|}$ is
the inverse of the exponential of a subordinator with no drift, in
particular it is a pure jump process. Furthermore it is clear by
construction that all the jumps of $B_i(.)$ correspond to some
atom of $\tilde{K}$

For each $t$ define $\psi_t(.)$ the application from $\N
\rightarrow \N$ inverse of $\phi(t,.)$, i.e.
$$\psi_t \left( \phi_t \left( i \right)\right)=i$$ (exists
because $\phi$ is surjective on $\N$).

\begin{lemma}\label{saut pur}
Under the above assumption on $\nu$,
\begin{itemize}
\item for all $k > 0$,  $\lambda_1(t)+ \lambda_2(t) +...+
\lambda_k(t)$ is a pure jump process.
\item with probability one, for all $t\geq 0$, if $t$ is an atom
for $\lambda$ then $\tilde{K}$ has an atom at t.
\end{itemize}
\end{lemma}

\begin{proof}

We will begin by proving the result for $\lambda_1$, the size of
the largest fragment and then turn our attention to the small
ones.

$\lambda_1$ is a supremum of a countable family of pure jump
processes (the $|B_i(.)|$). However it is easy to exhibit an
example of a supremum of a countable family of pure jump processes
that is not a pure jump process. So the proof will consist in
showing that almost surely on a fixed time interval $\lambda_1$ is
the supremum of a finite number of pure jump processes.

For this proof only, it is convenient to work with so-called
\textit{interval fragmentations}.

\textit{Interval fragmentations} are a particular case of
\textit{object fragmentations} that we presented in the
introduction for which the "object" $E$ is simply the interval
$[0,1]$ endowed with the Lebesgue measure. More precisely, call
$\mathbf{\nu}$ the space of the open subsets of $[0,1]$. Elements
of $\mathbf{\nu}$ admit a unique decomposition in intervals (in
the sense that the ordered vector of the lengths is unique). An
interval decomposition is a process $F(t)$ with values in
$\mathbf{\nu}$ such that for any $0 \leq s <t$ one has $F(s)
\subseteq F(t)$ i.e. $F(t)$ is finer than $F(s)$.

Take a sequence $(u_i)_{i \in \N}$ of iid variables uniformly
distributed on $[0,1]$. $F$ is then transformed into a
$\P$-process $\Pi$ by the following rule $$i
\overset{\Pi(t)}{\sim} j \Leftrightarrow [u_i,u_j]\subseteq
F(t).$$ This last process obviously conserves the refinment
property, moreover, if we define interval fragmentations to have a
scaling and branching property, $\Pi$ will be a
$\P$-fragmentation.

We refer to \cite{bertoin2} for a precise definition of interval
fragmentation and the equivalence between interval fragmentations
and partition fragmentations.

There is no loss of generality in supposing that $\Pi$ is
constructed from an interval fragmentation $F(t)$ and  a sequence
$(u_i)_{i \in \N}$ of iid variables uniformly distributed on
$[0,1]$.

Denote $(I_i(t) , i \in \N)$ the associated ordered length of the
interval decomposition of $F$ (which are also the associated
ordered frequencies of $\Pi(t)$ ). If $I^{(i)}(t)$ denote the
length of the interval that contains $u_i$ in the interval
decomposition of $F(t)$, then $$I^{(i)}(t) = l_i(t)$$ where
$l_i(t)=|\beta(i,t)|$ is the asymptotic frequency of the block of
$\Pi(t)$ that contains $i$.

Calling $\tau_n$ the stopping time $\inf\{t>0, |B_n(t)|>0\}$ we
have that at $\tau_n$ $$\forall i<n , n
\overset{\Pi(\tau_n)}{\not\sim} i,$$ thus $u_n$ does not belong
to any block of $F(\tau_n)$ that contains some $u_i$ for any $i
<n$, hence the asymptotic frequency of the block of $F(\tau_n)$
that contains $u_n$ is bounded from above by $\sup_{i,j \in
\{1,..,n\}}|u_i-u_j|$ which converge to $0$ almost surely when $n
\rightarrow \infty$.

Note that $$\sup_{r \in \R^+ } \{ |B_n(r)|\}=|B_n(\tau_n)|$$ to
see that $$\lim_{n \rightarrow \infty}\left( \sup_{r>0}\left(
|B_n(r)| \right) \right)=0 \text{ a.s.}$$

Now fix $\epsilon>0$ and $n_0$ and condition on the events $\{
\lambda_1(T) \geq \epsilon \}$, and $$\{ \sup_{n>n_0}\{
\sup_{r>0}\left( |B_n(r)| \right) \} \}<\epsilon \}.$$ Note that
the probability of the second event can be taken arbitrarily
close to 1 by taking $n_0$ sufficiently large. On this event, for
all $r \in [0,T]$ we have that
$$\lambda_1(r)=\max_{i=1,...n_0}|B_i(r)|.$$ Thus $\lambda_1(.)$
is a pure jump process because all the $|B_i(.)|$ are. Moreover
$\lambda_1(.)$ only jumps at times at which $\tilde{K}$ has an
atom for the same reason.

We now turn our attention to the other fragments.

Let $0<a<b$ and suppose the result is proven for the $(\lambda_i,
i\in \{1,...,k-1\})$. At time $a$ there is almost surely an
infinite number of blocks each with a positive asymptotic
frequency, suppose
\begin{eqnarray*}
\lambda_k(a)= |B_{\psi_a(k)}(a)|>\epsilon>0.
\end{eqnarray*}

Call a child of $B_{\psi_a(1)}(a)$ a block of $\Pi(a+u), u>0$
included in $B_{\psi_a(1)}(a)$. Denote by $C_1(b)$ (resp.
$C_j(b)$) the size of the largest child of $B_{\psi_a(1)}(a)$ at
time $b$ (resp. the size of the largest child of
$B_{\psi_a(j)}(a)$ at time $b$). They are almost surely strictly
positive. Let $\eta>0$ and condition on the event
$$C_1(b)\wedge C_2(b)\wedge.. \wedge C_k(b)>\eta$$ then it is clear that
$\eta$ is a lower bound for $\inf_{t \in [a,b]} \{\lambda_k(t)\}$,
thus the same argument as in the $\lambda_1$ case allow us to
consider only a finite number of fragment to be sure to "catch"
$\lambda_k$. More specifically, conditioned on the event $$\left\{
\sup_{n>n_0}\left\{ \sup_{u>0}\left\{\left( |B_n(u)| \right)
\right\} \right\}<\eta \right\},$$ whose probability can be
controlled through $n_0$ to be as close as we wish to 1, we can
write $\lambda_k$ as
$$\forall u \in [a,b] , \lambda_k(u)=\sup
\left(\{|B_j(u)|\}_{j=1,...,n_0} \backslash
\{\lambda_1(u),\lambda_2(u),...,\lambda_{k-1}(u)\}\right).$$ As
the $|B_i|$ and the $\lambda_1,...,\lambda_{k-1}$ are pure jump
processes, $\lambda_k$ is a pure jump process on $[a,b]$, and its
jumps correspond to atoms of $\tilde{K}$ for which
$\tilde{k}(t)\leq k$ (and these atoms are themselves images of
atoms of $K$ for which $k(t)\leq n_0$). So by induction the result
is proven.

\end{proof}

In conclusion, if we call $\Gamma$ the set of times at which
$(\Delta(t),k(t))$ has an atom. Then writing
$\lambda(t)=(\lambda_1(t),\lambda_2(t),...)$ for $\Lambda(\Pi(t))$
:
\begin{enumerate}
\item  $\lambda(.)$ is a pure jump process, c\`adl\`ag and starts almost surely from $(1,0,0,...)$
\item  if $t \not\in \Gamma$,
$$\lambda(t)=\lambda(t_{-})$$

\item if $t$ is a jump-time for
$\lambda$, then almost surely $t \in \Gamma$ and $\lambda(t)$ is
the reordering of the concatenation of two sequences :
$(\lambda_i(t_-))_{\{i \neq k(t)\}}$ and
$\lambda_{\phi_t(k(t))}(t_-)\Lambda(\Delta(t))$. \
\end{enumerate}

As $\lambda$ is a pure jump process it is completely defined by
this description.

All we have to do now is collect the preceding results : let
$K=(\Delta(t),k(t))$ be a Poisson point process with measure
intensity $\mu_{\nu} \otimes \#$ and let $\Pi$ the associated
$(\nu,0)$ homogeneous $\P$-fragmentation. Then the Poisson point
process $\left(\Lambda(\Delta(t)),\phi_{t_-}\left( k(t)
\right)\right)$ and the asymptotic frequency process
$\Lambda(\Pi(t))$ have the desired properties, so Theorem
\ref{prop saut pur} is proved.

\section{Small time Asymptotic behavior}

In this section we use the Poisson construction of ranked
fragmentations we just established to study their asymptotic
behavior near 0. The results we give are very close in spirit to
those concerning the asymptotic behavior of subordinators.

A subordinator, say $\xi$, is an increasing L\'evy process whose
distribution is specified by its Laplace exponent $\Psi$ that is
given by the identity
$$\mathbb{E}(\exp{\{-q\xi_t\}})=\exp{\{-t\Psi(q)\}}$$ and the
L\'evy-Khintchine formula $$\Psi(q)=k+dq+\int_{\left] 0,\infty
\right[} (1-e^{-qx})\upsilon(dx)$$ where $k\geq 0$ is the
so-called killing rate, $d\geq 0$ is the drift coefficient and
$\upsilon$ a measure on $\left] 0, \infty \right[$ with $\int(1
\wedge x) \upsilon(dx) < \infty$, called the L\'evy measure of
$\xi$.

The asymptotic behavior of these processes is well known, for
instance we have results concerning their distribution :
$$\frac{1}{t}P(\xi(t) \in .) \underset{t \rightarrow 0+}{\rightarrow} \upsilon(.)$$ (see Corollary 8.9 in \cite{sato}).

On the other hand, under conditions of regular variation on the
tail of $\upsilon$, there are also results concerning the sample
path behavior of the limsup and the law of the iterated logarithm
(see for instance the end of chapter III in \cite{bertoin5}). More
precisely :
\begin{itemize}
\item (law of the iterated logarithm) A necessary an sufficient condition for the Laplace
exponent $\Psi$ of $\xi$ to be  regularly varying near $\infty$
with index $a \in (0,1)$ is that the drift coefficient is 0 and
$\bar{\upsilon}(x)=\upsilon(\left] x,\infty \right[)$ is regularly
varying in $0+$ with index $-a$. In this case it holds with
probability 1 that
$$\liminf_{t\rightarrow 0+} (\frac{\xi(t) \Psi^{-1}(t^{-1}\log
|\log t|)}{\log |\log t|})=a(1-a)^{(1-a)/a}.$$
\item suppose the drift is 0 and let $h
:\left[0,\infty\right)\rightarrow \left[0,\infty\right)$ be an
increasing function such that the function $t\rightarrow h(t)/t$
increases as well. Then the following assertions are equivalent :
\begin{enumerate}
\item a.s. $$\limsup_{t \rightarrow 0+} (\xi(t)/h(t))=\infty ;$$
\item $$\int_0^1 \bar{\nu}(h(t))dt=\infty ;$$
\item $$\int_0^1 \left\{
\Psi(1/h(t))-(1/h(t))\Psi'(1/h(t))\right\}dt=\infty.$$
\end{enumerate}
Finally if these assertions fail to be true, then almost surely
$$\lim_{t \rightarrow 0+} (\xi(t)/h(t))=0.$$
\end{itemize}

Thus to study the asymptotic behavior of a fragmentation we may
benefit from the fact that $|B_1|$ (the mass of the block that
contains 1) can be described in terms of a subordinator (see
\cite{bertoin1}).

We focus on the behavior of the largest ($\lambda_1$) and of the
second block ($\lambda_2$) of a ranked fragmentation eventhought
we have more general result in the case of so-called
\textit{binary} fragmentations.

Although the study of $\lambda_1$ is relatively straightforward,
$\lambda_2$ requires to use some results of the record-processes
theory. Most of those that will be used in this section are well
known or are adapted from standard facts that can be found in
most textbooks on the matter. See \cite{bingham} for instance.

First note that $\lambda_2(t)$ is not monotone, more precisely it
decreases when the second largest fragment undergoes a dislocation
and can increase when the largest fragment undergoes a dislocation
and one of the new fragment created becomes the second largest.

The idea is to use the Poisson construction : near 0 the largest
fragment is almost of size 1, thus the second largest fragment is
always a "direct son" of the main one, and we shall be able to
express its law in terms of the distribution of the largest
fragment that has detached from the main.

For a general $\R$-valued P.P.P. $K=(K_t, t \geq 0)$ with
intensity measure $\mu$ such that  $\forall \epsilon
>0, \mu(\left]\epsilon,\infty\right])<\infty$, it is possible to
define the associated record process $R(t)$ as follows : at time
$t$
$$R(t)=\max_{s\leq t}\{K_s\}.$$

Let $\lambda$ be a homogeneous $\sfleche$ fragmentation with
characteristic $(\nu,c)$ constructed from the P.P.P.
$$K=(S(t),k(t))_{t\geq
0}=((s_1(t),s_2(t),...),k(t))_{t\geq 0}$$ of intensity measure
$\nu \otimes \#$. Let $(S^{(i)}(t), t \geq 0)= (s_j^{(i)}(t),
j=1,2,... ; t \geq 0)$ be  the P.P.P. with values in $\sfleche$
derived of $K$ by keeping the points such that $k(t)=i$ (the
second coordinate being always $i$, it is not expressed). So
$s^{(i)}_j(t)$ is the relative size of the $j^{th}$ block of the
dislocation occurring at time $t$ on the $i^{th}$ block. $S^{(i)}$
is a P.P.P. with intensity measure $\nu$. The $\R$-valued point
process $(s^{(i)}_j(t))$ is thus a P.P.P. with intensity
$$\nu_j(dx)=\nu(\{s=(s_1,s_2,...) \in \sfleche : s_j \in dx
\}).$$

Introduce the function
$$x \rightarrow \overline{\nu}_2(x)=\nu(s \in \sfleche : s_2 \geq
x)$$ from $[0,\frac{1}{2}] \rightarrow \R^+$, and denote by $f$
its generalized inverse.

Note that $\overline{\nu}_2(.)$ is finite, i.e. for all $x>0$
$\overline{\nu}_2(x)<\infty$. To see this, let $b \in [0,1/2]$
\begin{eqnarray*}
\int_{\sfleche}(1-s_1)\nu(ds) &\geq& \int_{\sfleche}s_2\nu(ds)\\
 &=& \int_0^{1/2}x\nu_2(dx) \\ &\geq& \int_{b}^{1/2}x\nu_2(dx)
 \\ &>& b\overline{\nu}_2(b)
 \end{eqnarray*}

Let $R(t)$ designate the record at time $t$ of the P.P.P.
$s^{(1)}_2(.)$ which is well defined according to the above
argument.

\begin{proposition}\label{small time behavior}
Let $$\lambda=(\lambda(t),t\geq
0)=(\lambda_1(t),\lambda_2(t),...),t\geq 0)$$ be a homogeneous
$\sfleche$ fragmentation with characteristic $(\nu,c)$, then
\begin{enumerate}
\item there exists a subordinator $\xi$ with drift $c$ and L\'evy
measure $$L(dx)=e^{-x}  \nu(-\log s_1 \in dx) , x \in
\left]0,\infty \right[$$ such that  $$1-\lambda_1(t) = 1 -\exp{
\xi(t)} $$ for $t$ small enough a.s.
\item $$\lambda_2(t) \sim R(t), \qquad t \rightarrow 0+ \qquad  a.s.$$
\end{enumerate}
\end{proposition}

\begin{proof}\textit{(Proposition \ref{small time behavior}-(1))} Assume that
$c=0$, then consider $\tilde{\nu}$ the image of $\nu$ by the
application $\sfleche \rightarrow \sfleche : (s_1,s_2,...)
\rightarrow (s_1,0,...)$. Let $\tilde{\lambda}$ a homogeneous
$(\tilde{\nu},0)$ $\sfleche$-fragmentation, which thus has no
erosion and almost surely for all $t$ only one block that has
positive mass. There is no loss of generality in supposing
$\tilde{\lambda}=\Lambda(\tilde{\Pi})$ where
$\tilde{\Pi}=(\tilde{B_1},\tilde{B_2},...)$ is a homogeneous
$(\tilde{\nu},0)$ $\P$-fragmentation. Define
$\tilde{\xi}(t)=-\log(\tilde{B_1}(t))$. As long as
$\tilde{\xi}(t)<\infty$ we have that
$\tilde{\xi(t)}=-\log(\tilde{\lambda_1}(t))$ (because it is the
only block which is not reduced to a singleton).

Next we condition on the event $\lambda_1(t) \geq 1/2$, for any $s
\leq t , \lambda_1(s)$ is either $\lambda_1(s_-)$,or the largest
fragment issued from a dislocation of $\lambda_1(s_-)$. By right
continuity $P(\lambda_1(t) \geq 1/2) \convtenzero 1$. On
$\lambda_1(t) \geq 1/2$ one has that $\forall s \leq t ,
\lambda_1(s)=\tilde{\lambda_1}(s)$ (one can construct $\lambda$
and $\tilde{\lambda}$ from $K$ (the same PPP) using its image by
the above mentioned transformation for $\tilde{\lambda}$). Thus
conditionally on $s \leq t$
$$1-\lambda_1(s)=1-\exp{(-\tilde{\xi}(s))}.$$
\end{proof}

This equivalence relation combined with subordinator properties
have immediate consequences such as
$$\frac{1}{t}P(1-\lambda_1(t)>x) \convtenzero L([-\log(1-x),1])$$
where $L(dx)=e^{-x}\nu(-\log(s_1) \in dx).$

For the second point the idea is to describe the asymptotic
behavior of $\lambda_2$ in terms of the records of $s^{(1)}_2$. We
begin with the following technical lemma
\begin{lemma}
\label{lem1} Let
$$\chi_t=\left(\prod_{u\in \left[0,t\right[}s^{(1)}_1(u)\right)
\left(\prod_{u\in \left[0,t\right[}s^{(2)}_1(u)\right), $$ and
suppose $c=0$ (there is no erosion) then on the event $\{
\lambda_1(t) \geq 1/2 \}$
$$\chi_t R(t) \leq  \lambda_2(t) \leq R(t).$$
\end{lemma}
\begin{proof}
As noted earlier, one can suppose that $\lambda(.)$ is the
asymptotic frequency of some $(\nu,0)$ $\P$-fragmentation $\Pi$,
and $ K$ is the image of the P.P.P. $$(\Delta(.),k(.))
\rightarrow (\Lambda(\Delta(.)),\phi(.,k(.)))$$ with intensity
measure $(\mu_{\nu} \otimes\# )$. At time $t$ we recall the
notation $\psi_t(1)$ for the least element of the block of
greatest asymptotic frequency in $\Pi$, which is well defined.

Fix $t$, and consider $(\Pi(t-u))_{u \in [0,t]}$, the
fragmentation where the time have been reversed. Informally it is
a coalescence, whose final state at $u=t$ is almost surely the
trivial partition and which is left continuous. Thus the
functional $$f_i(u)=\mathbf{1}_{\{i \overset{\Pi(t-u)}{\sim}
\psi_t(1)\}}$$ (that is $f_i(u)$ is $1$ if at time $(t-u)$ if $i$
is in the same block that the integer which is the least element
of the largest block at time $t$ and $0$ otherwise) is left
continuous and $f_i(t)=1$ a.s. Thus, almost surely
$$D_i(t)=t-\sup \{ u\in [0,t]: f_i(u)=0 \} < t.$$ Note that as we
are on $\{\lambda_1(t) \geq \frac{1}{2}\}$, $\psi_t(1)$ is always
in the block of greatest asymptotic frequency of $\Pi(t-u)$ for
any $ u \in [0,t]$; so $D_i(t)$ is the \textit{detachment time} of
$i$ from the main block (if $i$ is still in the main block,
$D_i(t)$ is taken equal to $t$).

Now take $k \geq 2$, and suppose that at time $t$ there is at
least $k$ blocks (almost surely the case under our hypothesis) so
$\psi_t(k)$ (the least element of the block of $k$-th greatest
asymptotic frequency in $\Pi$ at time $t$) is well defined,
almost surely
$$D_{\psi_t(k)}(t)>0.$$ so if we note $\beta(i,u)$ for the block of
$\Pi(u)$ that contains $i$ and $D(k,t)=D_{\psi_t(k)}(t)$ we have
that

\begin{eqnarray*}
\beta((\psi_t(k)),D(k,t)_-)&=&
\beta(\Psi_{D(k,t)_-}(1),D(k,t)_-) \\
|\beta\left(\left(\psi_t(k)\right),D(k,t)\right)|&<&
|\beta(\Psi_{D(k,t)}(1),D(k,t))|
\end{eqnarray*}
(recall that $\beta(\Psi_t(1),t)$ is the largest block at time
$t$). Thus $\lambda_k(t) \leq s^{(1)}_2(D_{\psi_t(k)}(t))$. As
obviously $$ s^{(1)}_2(D_{\psi_t(k)}(t)) \leq R(t)$$ we conclude
that
\begin{eqnarray}
\label{majo} \lambda_k(t)  \leq  R(t).
\end{eqnarray}

We now prove the lower-bound part of the lemma.

Let $T(t)=\inf\{u\leq t : R(u)=R(t)\}$ (the "record-time"). Note
that for all $u\in [0,t]$ at which $S^{(2)}$ has an atom,
$$\lambda_{2}(u_-)s^{(2)}_1(u)\leq \lambda_2(u),$$ this is not an
equality because the largest fragment issued of the dislocation of
$\lambda_2(u_-)$ can be smaller than $\lambda_3(u_-)$. Then, for
all $u\in [0,t]$ not an atom for $S^{(2)}$,$$\lambda_{2}(u-)\leq
\lambda_2(u),$$ this is due to the fact that $u$ could be an atom
for $S^{(1)}$, for which
$\lambda_1(u_-)s^{(1)}_2(u)>\lambda_2(u_-)$. Recalling that we are
still conditioning on $\{ \lambda_1(t) > \frac{1}{2}\}$ we have,
using the fact that $\lambda_2$ is a pure jump process, that
\begin{eqnarray}
\label{majo2} \lambda_{2}(T(t)) \bigg(\prod_{u \in \left[
T(t),t\right[} s_1^{(2)}(u) \bigg)\leq \lambda_2(t)
\end{eqnarray} and here again this is not an equality
because a reordering might occurs.

Then remark
\begin{eqnarray}
\label{remark} \lambda_{2}(T(t)) \geq R(t)\bigg( \prod_{u \in
\left[0, T(t)\right[} s_1^{(1)}(u) \bigg) .
\end{eqnarray}
In words : at the time of the record $R(1,t)$, the second
fragment issued of the dislocation of $\lambda_1$, is not
necessarily $\lambda_2$, but in any case it is smaller or equal.

We can combine (\ref{majo2}) and (\ref{remark}) to get
\begin{eqnarray}
\label{mino} \chi_t R(t) \leq \lambda_t(2)
\end{eqnarray}

\end{proof}

We can now prove the second part of proposition \ref{small time
behavior} :

\begin{proof}\textit{Proposition \ref{small time behavior}-(2)} When
$c=0$ we now only have to show that $\chi_t \convtenzero 1$ almost
surely. $\bigg(\prod_{u \in \left[0,
t\right[}(s^{(1)}_1(u))\bigg)$ and $\bigg(\prod_{u \in \left[0,
t\right[}(s^{(2)}_1(u))\bigg)$  are independent and identically
distributed, and on the event $\bigg(\prod_{u \in \left[0,
t\right[}(s^{(1)}_1(u))\bigg) \geq \frac{1}{2}$ this last quantity
is exactly the $\lambda_1(t)$ of some $(\nu,0)$ fragmentation,
thus almost surely $$\bigg(\prod_{u \in \left[0,
t\right[}(s^{(1)}_1(u))\bigg) \convtenzero 1$$ which thus
concludes our proof.

Finally in the case of a homogeneous $(\nu,c)$ fragmentation
$\lambda$ with $c \geq 0$, the effect of the erosion is just of
multiplying the size of each fragment by a factor $e^{-ct}$. So
clearly the upper bound of Lemma \ref{lem1} is still valid, on
the other hand we have
$$\xi(t)e^{-ct}R(t) \leq \lambda_2(t)$$ and only a slight
modification of the proof for the case $c=0$ is needed.
\end{proof}

\begin{remark}
If $\overline{\nu}_2(.)$ is regularly varying with index $(-a)$ in
$0^+$, $a \geq 0$, classical results of record-processes theory
used with proposition \ref{small time behavior} show that
$$\frac{\lambda_2(t)}{f(\frac{1}{t})}
\overset{\mathcal{L}}{\rightarrow}L$$ when $t \searrow 0$ where
$L$ is the extreme law with distribution function $\exp(-x^{-a}).$
\end{remark}

\begin{remark}
Suppose that $\lambda$ is a binary fragmentation, that is $\nu$
has its support in the subset of $\sfleche$ defined as $\{s \in
\sfleche , s_3=s_4=...=0\}$  and  that
$\overline{\nu}_2(x)=\nu(\{s \in \sfleche : s_2 \geq x\})$ is
regularly varying near 0 with index $-a$. Then using the same
ideas as in the above arguments one can show that we have the
following asymptotic distributions of the renormalized
$\lambda_k$ for any $k>1$  : $$\forall k >1 , a.s. \qquad
\lambda_k(t)\underset{t \searrow 0^+}{\sim}R_2(k,t).$$ As a
consequence $$\frac{\lambda_k(t)}{f(1/t)}
\overset{\mathcal{L}}{\rightarrow} L(k,a)$$ where $L(k,a)$ is the
law with repartition function $$F_{k,a}(x)=\left( \sum_{i\in
\left[0, k-1\right]}e^{-x^{-a}}\frac{(x^{-ai})}{i!} \right)$$ and
$f$ is the generalized inverse of $x \rightarrow \bar{\nu}_2(x)$.
More generally, the convergence in law holds jointly, the limit
distribution function for the $N$ largest blocks being given by
$$ f_N(x_2,x_3,...x_N) = \left(
\prod_{i=2}^{i=N-1}\left(\exp{-x_i^{-a}}\right)\right)\int_0^{x_N}(\exp{-u^{-a}})
\nu (s_2 \in du )   $$ for $x_1>x_2>...>x_N$ (see \cite{perman}
for instance).

\end{remark}

\begin{remark}
In the case where the fragmentations considered are not
homogeneous but only self-similar and without erosion, a slightly
more technical version of Theorem \ref{prop saut pur} still
stands : i.e. it is possible to give an explicit Poisson
construction of any $(\alpha,\nu,0)$ $\sfleche$-fragmentation.
This allow us to extend the results of section 4 to the case of a
self-similar fragmentation with index $\alpha>0$.
\end{remark}

\end{document}